\documentclass{amsart}
\usepackage{amssymb}
\usepackage{amsfonts}
\usepackage{latexsym}

\newtheorem{theorem}{Theorem}[section]
\newtheorem{lemma}[theorem]{Lemma}
\newtheorem{proposition}[theorem]{Proposition}
\newtheorem{corollary}[theorem]{Corollary}
\theoremstyle{definition}

\newtheorem{example}[theorem]{Example}

\newtheorem{remark}[theorem]{Remark}

\newcommand{\id}{\text{id}}

\newcommand{\Fun}{\text{Fun}}
\newcommand{\Irr}{\text{\rm Irr}}
\newcommand{\FPdim}{\text{\rm FPdim}}
\newcommand{\Ext}{\text{\rm Ext}}

\newcommand{\gr}{\text{gr}}

\newcommand{\Rep}{\text{Rep}}

\newcommand{\ot}{\otimes}

\newcommand{\ben}{\begin{enumerate}}
\newcommand{\een}{\end{enumerate}}

\newcommand{\Rad}{{\text{Rad}}}

\begin{document}

\title{Finite dimensional quasi-Hopf algebras with radical of
codimension 2}
\author{Pavel Etingof}
\address{Department of Mathematics, Massachusetts Institute of Technology,
Cambridge, MA 02139, USA} \email{etingof@math.mit.edu}

\author{Shlomo Gelaki}
\address{Department of Mathematics, Technion-Israel Institute of
Technology, Haifa 32000, Israel}
\email{gelaki@math.technion.ac.il}

\maketitle

\section{Introduction}

It is shown in \cite{eo}, Proposition 2.17, that a finite
dimensional quasi-Hopf algebra with radical of codimension 1 is
semisimple and 1-dimensional. On the other hand, there exist
quasi-Hopf (in fact, Hopf) algebras, whose radical has codimension
2. Namely, it is known \cite{n} that these are exactly the Nichols
Hopf algebras $H_{2^n}$ of dimension $2^n$, $n\ge 1$ (one for each
value of $n$).

The main result of this paper is that if $H$ is a finite
dimensional {\bf quasi-Hopf algebra} over $\mathbb{C}$ with
radical of codimension 2, then $H$ is twist equivalent to a
Nichols Hopf algebra $H_{2^n}$, $n\ge 1$, or to a lifting of one of
the four special quasi-Hopf algebras $H(2)$, $H_+(8)$, $H_-(8)$,
$H(32)$ of dimensions 2, 8, 8, and 32, defined in Section 3. As a
corollary we obtain that any finite tensor category which has two
invertible objects and no other simple object is equivalent to
$\Rep(H_{2^n})$ for a unique $n\ge 1$, or to a deformation of the
representation category of $H(2)$, $H_+(8)$, $H_-(8)$, or $H(32)$.
In the case of $H(2)$ such a lifting (deformation) must clearly be
trivial; for the other three cases we plan to study possible
liftings in a later publication.

As another corollary we prove that any nonsemisimple quasi-Hopf
algebra of dimension 4 is twist equivalent to $H_4$.

Thus, this paper should be regarded as a beginning of the
structure theory of finite dimensional basic quasi-Hopf algebras.
It is our expectation that this theory will be a nontrivial
generalization of the deep and beautiful theory of finite
dimensional basic (or, equivalently, pointed) Hopf algebras (see
\cite{as} and references therein).

{\bf Acknowledgments.} We thank V. Ostrik for suggesting to
generalize the result of Nichols to the quasi-Hopf case.

The first author was partially supported by the NSF grant
DMS-9988796. The second author's research was supported by
Technion V.P.R. Fund - Dent Charitable Trust - Non Military
Research Fund, and by The Israel Science Foundation (grant No.
70/02-1). He is also grateful to MIT for its warm hospitality.
Both authors were partially supported by the BSF grant No.
2002040.

\section{Preliminaries}

All constructions in this paper are done over the
field of complex numbers $\mathbb{C}$.

\subsection{} Recall that a finite rigid tensor category ${\mathcal C}$ over
$\mathbb{C}$ is a rigid tensor category over $\mathbb{C}$ with
finitely many simple objects and enough projectives such that the
neutral object ${\bf 1}$ is simple (see \cite{eo}). Let
$\Irr({\mathcal C})$ denote the (finite) set of isomorphism
classes of simple objects of ${\mathcal C}$. Then to each object
$X\in {\mathcal C}$ there is attached a positive number
$\FPdim(X)$, called the Frobenius-Perron (FP) dimension of $X$,
and the Frobenius-Perron (FP) dimension of ${\mathcal C}$ is
defined by
$$\FPdim({\mathcal C})=\sum_{X\in \Irr({\mathcal C})}
\FPdim(X)\FPdim (P_X),$$ where $P_X$ denotes the projective cover
of $X$ \cite{e}. By Proposition 2.6 in \cite{eo}, ${\mathcal C}$
is equivalent to $\Rep(H)$, $H$ a finite dimensional quasi-Hopf
algebra, if and only if the FP-dimensions of its objects are
integers. We refer the reader to \cite{d} for the definition of a
quasi-Hopf algebra and the notion of twist (= gauge) equivalence
between quasi-Hopf algebras.

\subsection{} Recall that Nichols' Hopf algebra $H_{2^n}$
($n\ge 1$) is the $2^n-$dimensional Hopf algebra generated by
$g,x_1,\dots,x_{n-1}$ with relations $g^2=1$, $x_i^2=0$ ($1\le
i\le {n-1}$), $x_ix_j=-x_jx_i$ ($1\le i\ne j\le n-1$) and
$gx_i=-x_ig$, where $g$ is a grouplike element and
$\Delta(x_i)=x_i\ot 1 + g\ot x_i$ \cite{n}. The Hopf algebra $H_4$
is known as the Sweedler's Hopf algebra \cite{s}. Nichols' Hopf
algebras are the unique (up to isomorphism) Hopf algebras with
radical of codimension $2$. It is well known (and easy to check)
that $H_{2^n}$ is self dual.

\subsection{} The following is the simplest (and well known) example of
a quasi-Hopf algebra not twist equivalent to a Hopf algebra. The
2-dimensional quasi-Hopf algebra $H(2)$ is generated by a
grouplike element $g$ such that $g^2=1$, with associator
$\Phi=1-2p_-\otimes p_-\otimes p_-$ (where $p_+:=(1+ g)/2,
p_-:=(1-g)/2$), distinguished elements $\alpha=g$, $\beta=1$, and
antipode $S(g)=g$. It is not twist equivalent to a Hopf algebra,
and any 2-dimensional quasi-Hopf algebra is well known to be twist
equivalent either to $\Bbb C[\Bbb Z_2]$ or to $H(2)$.

\subsection{} Let $H$ be a finite dimensional algebra,
and let $I:=\Rad(H)$ be the Jacobson radical of $H$. Then we have
a filtration on $H$ by powers of $I$. So one can consider the
associated graded algebra $\gr(H)=\bigoplus_{k\ge 0} H_k$,
$H_k:=I^k/I^{k+1}$ ($I_0:=H$).

The following lemma is standard.

\begin{lemma}\label{sta} (i) Let $L$ be any linear complement
of $I^2$ in $H$. Then $L$ generates $H$ as an algebra.

(ii) $\gr(H)$ is generated by $H_0$ and $H_1$.
\end{lemma}

\begin{proof} (i) Let $B\subseteq H$ be
the subalgebra generated by $L$. We will show by induction in $k$
that $B+I^k=H$ for all $k$. This implies the statement, since
$I^N=0$ for some $N$. The base of induction, $k=2$, is clear. So
assume that $k>2$. By the induction assumption, all we need to
show is that $I^{k-1}\subseteq B+I^k$. So let us take $a\in
I^{k-1}$, and show that $a\in B+I^k$. We may assume that
$a=i_1\cdots i_{k-1}$, where $i_m\in I$ for all $m$. Let $\bar
i_m\in L\cap I$ be the unique elements such that $\bar i_m-i_m\in
I^2$. Then $\bar a:=\bar i_1\cdots \bar i_{k-1}\in B$, and $a-\bar
a\in I^k$. We are done.

(ii) Clearly, $\gr(I)$ is the radical of $\gr(H)$. Thus we only
need to show that $\gr(I^2)=\gr(I)^2$ (then everything follows
from (i)). The inclusion $\gr(I^2)\supseteq \gr(I)^2$ is clear, so
let us prove the opposite inclusion. Since both spaces are graded,
it suffices to pick a homogeneous element $a\in \gr(I^2)$ and show
that it lies in $\gr(I)^2$. Suppose $a$ has degree $k$, i.e. $a\in
I^k/I^{k+1}$. Let $\tilde a$ be a lifting of $a$ to $I^k$. Then
there exist elements $i_{1j},...,i_{kj}\in I$, $j=1,...,m$, such
that $\tilde a=\sum_j i_{1j}\cdots i_{kj}$. Let $\bar
i_{1j},...,\bar i_{kj}$ be the projections of $i_{1j},...,i_{kj}$
to $I/I^2\subseteq \gr(I)$. Then $a=\sum_j \bar i_{1j}\cdots \bar
i_{kj}$. Since $k\ge 2$, we are done.
\end{proof}

\subsection{}
Let $H$ be a basic quasi-Hopf algebra (i.e., all irreducible
$H$-modules are $1$-dimensional) and let $I:=\Rad(H)$ be the
radical of $H$. Then $I$ is a quasi-Hopf ideal (i.e.
$\Delta(I)\subseteq H\ot I + I\ot H$, $S(I)=I$ and
$\varepsilon(I)=0$), so we have a quasi-Hopf algebra filtration on
$H$ by powers of $I$. The associated graded algebra
$\gr(H)=\bigoplus_k H_k$ is thus also a quasi-Hopf algebra. We
have $H_0=\Fun(G)$, where $G$ is a finite group (of characters of
$H$), $H_1=\bigoplus_{a,b\in G}\Ext^1(a,b)^*$. Note that by
definition, the associator of $\gr(H)$ lives in $H_0^{\otimes 3}$,
so it corresponds to an element of $H^3(G,\mathbb{C}^*)$.

\section{The main results}

Our first main theorem is the following.

\begin{theorem}\label{main1} Let $H$ be a finite
dimensional quasi-Hopf algebra over $\mathbb{C}$ with radical of
codimension $2$. Suppose that the associator of $H$ is trivial on
$1-$dimensional $H$-modules (i.e., the subcategory of $\Rep(H)$
generated by $1-$dimensional objects is equivalent to
$\Rep(\mathbb{Z}_2)$ with the trivial associator). Then $H$ is
twist equivalent to a Nichols Hopf algebra $H_{2^n}$, $n\ge 1$.
\end{theorem}

The proof of the theorem is given in Section 4.

The assumption on the triviality of the associator of $H$
on 1-dimensional modules is
clearly essential already for $\dim(H)=2$, because of the existence
of the quasi-Hopf algebra $H(2)$.
However, this assumption is also essential in higher dimensions,
as the following propositions show.

\begin{proposition}\label{prop8}
There exist two $8-$dimensional quasi-Hopf algebras $H_\pm (8)$
(permuted by the action of the Galois group), which have the
following structure. As algebras $H_\pm (8)$ are generated by
$g,x$ with the relations $gx=-xg$, $g^2=1$, $x^4=0$. The element
$g$ is grouplike, while the coproduct of $x$ is given by the
formula $$\Delta(x)=x\otimes (p_+ \pm ip_-) +1\otimes p_+x +
g\otimes p_-x,$$ where $p_+:=(1+ g)/2, p_-:=(1-g)/2$. The
associator is $\Phi=1-2p_-\otimes p_-\otimes p_-$, the
distinguished elements are $\alpha=g$, $\beta=1$, and the antipode
is $S(g)=g$, $S(x)=-x(p_+ \pm ip_-)$.
\end{proposition}

\begin{proposition}\label{prop32}
There exists a $32-$ dimensional quasi-Hopf algebra $H(32)$, which
has the following structure. As an algebra $H(32)$ is generated by
$g,x,y$ with the relations $gx=-xg$, $gy=-yg$, $g^2=1$, $x^4=0$,
$y^4=0$, $xy+iyx=0$. The element $g$ is grouplike, while the
coproducts of $x,y$ are given by the formulas $$\Delta(x)=x\otimes
(p_+ + ip_-) +1\otimes p_+x + g\otimes p_-x,$$
$$\Delta(y)=y\otimes (p_+ - ip_-) +1\otimes p_+y + g\otimes
p_-y,$$ where $p_+:=(1+ g)/2, p_-:=(1-g)/2$. The associator is
$\Phi=1-2p_-\otimes p_-\otimes p_-$, the distinguished elements
are $\alpha=g$, $\beta=1$, and the antipode is $S(g)=g$,
$S(x)=-x(p_+ + ip_-)$, $S(y)=y(p_+-ip_-)$. Thus, $H(32)$ is
generated by its quasi-Hopf subalgebras $H_+(8)$ and $H_-(8)$
generated by $g,x$ and $g,y$, respectively.
\end{proposition}

Propositions \ref{prop8} and \ref{prop32} are proved in Section 5.

Nonetheless, it turns out that outside of dimensions 2, 8 and 32
the assumption of the triviality of the associator of $H$ on
1-dimensional modules is satisfied automatically. This striking
fact follows from our second main theorem, which is the following.

\begin{theorem}\label{main3} Let $H$ be a finite
dimensional quasi-Hopf algebra over $\mathbb{C}$ with radical of
codimension $2$. Suppose that the associator of $H$ is nontrivial
on $1-$dimensional $H$-modules (i.e., the subcategory of $\Rep(H)$
generated by $1-$dimensional objects is equivalent to
$\Rep(H(2))$). Then $\gr(H)$ is twist equivalent to $H(2)$,
$H_+(8)$, $H_-(8)$, or $H(32)$; in particular, the dimension of
$H$ is $2$, $8$ or $32$. The quasi-Hopf algebras $H_+(8)$ and
$H_-(8)$ are not twist equivalent.
\end{theorem}

Theorem \ref{main3} is proved in Section 6.

\begin{corollary} Any finite tensor category which has two invertible
objects and no other simple object is tensor equivalent to
$\Rep(H_{2^n})$ for a unique $n\ge 1$, or to a deformation of the
representation category of $H(2)$, $H_+(8)$, $H_-(8)$, or $H(32)$.
\end{corollary}

\begin{proof} The FP-dimension of any object in $\mathcal{C}$
is an integer. Therefore, by Proposition 2.6 in \cite{eo},
$\mathcal{C}$ is equivalent to a representation category of a
finite dimensional quasi-Hopf algebra with radical of codimension
2. The result follows now from Theorems \ref{main1} and
\ref{main3}.
\end{proof}

\begin{corollary}\label{cor2} Any nonsemisimple quasi-Hopf algebra
of dimension 4 is twist equivalent to $H_4$.
\end{corollary}

\begin{proof} By dimension counting it follows that $H$ is basic;
that is, all irreducible representations of $H$ are
$1-$dimensional. Moreover, by Theorem 2.17 in \cite{eo}, $H$ has
more than one irreducible representation. Let $\chi$ be a
non-trivial $1-$dimensional representation of $H$, and let
$P_{\chi}$ denote the projective cover of $\chi$. Then
$P_{\chi}=\chi\ot P_{\varepsilon}$, where $P_{\varepsilon}$ is the
projective cover of the trivial representation $\varepsilon$.
Therefore $4$ equals $\dim(P_{\varepsilon})$ times the number of
simple objects in $\Rep(H)$. Since both numbers are greater than
$1$, we have that they are equal to $2$. So $H$ has $2$
irreducible representations $\varepsilon$ and $\chi$, with
$\chi^2=\varepsilon$. The result follows now from Theorems
\ref{main1} and \ref{main3}.
\end{proof}

\begin{corollary}\label{cor3}
A nonsemisimple finite tensor category ${\mathcal C}$ of
$\text{FP}$-dimension $4$ is tensor equivalent to
$\text{Rep}(H_4)$.
\end{corollary}

\begin{proof} It suffices to show that ${\mathcal C}$ has
integer FP-dimensions of objects; in this case by Proposition 2.6
in \cite{eo}, ${\mathcal C}$ is a representation category of a
4-dimensional quasi-Hopf algebra, and Corollary \ref{cor2}
applies.

Let $X$ be a simple object of ${\mathcal C}$ of FP-dimension $d$.
If $X\ot X^*\ne \bf 1$ then $X\ot X^*$ contains as constituents
$\bf 1$ and another object. So if $d\ne 1$ then $d\ge \sqrt{2}$.
On the other hand, the projective cover $P_{\bold 1}$ of the
neutral object has to involve other objects, so it is at least
2-dimensional. This shows that the only chance for ${\mathcal C}$
to have dimension 4 is when all simple objects are 1-dimensional.
We are done.
\end{proof}

\begin{remark}
We note that {\it semisimple} finite tensor categories of
FP-dimension 4 are easy to classify. Such a category is either the
category of modules over the group algebra of a group of order 4
with associativity defined by a 3-cocycle on this group, or a
Tambara-Yamagami category \cite{ty} corresponding to the Ising
model (see \cite{eno}, Proposition 8.32).
\end{remark}

\section{Proof of Theorem \ref{main1}}

By assumption, the associator of $\gr(H)$ is trivial. Therefore
$\gr(H)$ is twist equivalent to a Hopf algebra, and by the result
of Nichols \cite{n}, $\gr(H)$ is twist equivalent to a Nichols
Hopf algebra $A:=H_{2^n}$.

From this point we will identify $\gr(H)$ and $A$ as quasi-Hopf
algebras.

We will now show that $H$ itself is twist equivalent to a Nichols
Hopf algebra, completing the proof of the theorem.

Let $I_r$ be the radical of $H^{\otimes r}$. Then of course
$I_r=\sum_{k=1}^r H^{\otimes k-1}\otimes I\otimes H^{\otimes
r-k}$.

Let $\Phi$ be the associator of $H$. Then $\Phi=1+\phi$, where
$\phi\in I_3$. Assume that $\phi \in I_3^m$, but $\phi\notin
I_3^{m+1}$. We will show that by twisting we can change $\phi$ so
that it will belong to $I_3^{m+1}$. Then by a chain of twists we
can make sure that $\phi=0$, and by Nichols' result we are done.

Let $\phi'$ be the projection of $\phi$ to
$I_3^m/I_3^{m+1}=A^{\otimes 3}[m]$. Obviously, $\phi'$ is a
$3-$cocycle of $A^*$ with trivial coefficients. But
$A^*=A=\mathbb{C}[\mathbb{Z}_2]\ltimes \Lambda V$, and
$$H^3(A,\mathbb{C})= H^3(\mathbb{C}[\mathbb{Z}_2]\ltimes \Lambda
V,\mathbb{C})=H^3(\Lambda V,\mathbb{C})^{\mathbb{Z}_2}=
(S^3V^*)^{\mathbb{Z}_2}=0,$$ since $\mathbb{Z}_2$ acts on $V$ by sign.
Thus there exists an element $j'\in A^{\otimes 2}[m]$ such that
$dj'=\phi'$. Let $j$ be a lifting of $j'$ to $I_2^m$. Let
$J:=1+j$. Let $\Phi^J$ be the result of twisting $\Phi$ by $J$,
and $\Phi^J=1+\phi_J$. Then $\phi_J\in I_3^{m+1}$, as desired.

\section{Proof of Propositions \ref{prop8},\ref{prop32}}

\subsection{Proof of Proposition \ref{prop8}}

Let us first show that there exist quasi-Hopf algebras $H_\pm(\infty)$
defined in the same way as $H_\pm(8)$
but without the relation $x^4=0$.
To show this, it is sufficient to check
that $\Phi (\Delta\otimes \id)\Delta(x)=(\id\otimes
\Delta)\Delta(x)\Phi$ (the rest of the relations are
straightforward). This relation is checked by a simple direct
computation.

Now we must show that the ideal generated by $x^4$ in the
quasi-Hopf algebras $H_\pm(\infty)$ is a quasi-Hopf ideal.
To show this, we compute:
$$
\Delta(x^2)= x^2\otimes g+ [(1+i)(p_+\otimes p_++p_-\otimes
p_+)+(1-i)(p_+\otimes p_--p_-\otimes p_-)] (x\otimes x)+g\otimes
x^2,
$$
and hence
$$\Delta(x^4)=x^4\otimes 1+1\otimes x^4.
$$
So $x^4$ generates a quasi-Hopf ideal,
and the quotients $H_\pm(8):=H_\pm(\infty)/(x^4)$
(obviously of dimension 8) are quasi-Hopf algebras.
We are done.

\subsection{Proof of Proposition \ref{prop32}}

Let $H_{+-}(\infty)$ be the amalgamated product \linebreak
$H_+(8)*_{H(2)}H_-(8)$ (i.e. the algebra defined as $H(32)$ but
without the relation $xy+iyx=0$). This is obviously a quasi-Hopf
algebra (of infinite dimension). We must show that the principal
ideal in $H_{+-}(\infty)$ generated by $z:=xy+iyx$ is a quasi-Hopf
ideal. This follows from the easily established relation
$\Delta(z)=z\otimes 1+g\otimes z$.

Thus, $H(32):=H_{+-}(\infty)/(z)$ is a quasi-Hopf algebra. It is
easy to show that the elements $g^jx^ky^l$, $j=0,1$,
$k,l=0,1,2,3$, form a basis in $H(32)$. Thus, $H(32)$ is
32-dimensional.

\section{Proof of Theorem \ref{main3}}

Let $p_+,p_-$ be the primitive idempotents of $H_0$, and let
$g:=p_+-p_-$. Let $\Phi_0$ be the associator of $\gr(H)$. By using
a twist, we may assume that $\Phi_0=1-2p_-\otimes p_-\otimes p_-$
(the only, up to equivalence, nontrivial associator for
$\mathbb{C}[\mathbb{Z}_2]$).

Let $x$ be an element of $H_1$.
By Theorem 2.17 in \cite{eo},
$\Ext^1(\varepsilon,\varepsilon)=\Ext^1(\chi,\chi)=0$, thus
by subsection 2.5,
$$gxg^{-1}=-x.$$

Let $\chi$ be the
nontrivial character of $H_0$, and for $z\in H$ define
$$\xi(z):=(\chi\otimes
\id)(\Delta(z))\;\text{and}\;\eta(z):=(\id\otimes
\chi)(\Delta(z)).$$ Then $\xi(g)=\eta(g)=-g$.

For $x\in H_1$ we have $$\Delta(x)=x\otimes p_+ + p_+\otimes x + \eta(x)\otimes
p_- + p_-\otimes \xi(x).$$

By the quasi-coassociativity axiom, we have $\Phi_0 (\Delta\otimes
\id)\Delta(x)=(\id\otimes \Delta)\Delta(x)\Phi_0$.

\begin{lemma}\label{lem1}
The equation $\Phi_0 (\Delta\otimes \id)\Delta(x)=(\id\otimes
\Delta)\Delta(x)\Phi_0$ is equivalent to the relations
$$\xi^2(x)=-x,\;\eta^2(x)=-x,\;\text{and}\;\xi\eta(x)=
-\eta\xi(x).$$
\end{lemma}

\begin{proof}
For degree reasons the coassociativity equation is equivalent to the
system of three equations obtained by applying $\chi\otimes
\chi\otimes \id$, $\chi\otimes \id\otimes \chi$, and $\id\otimes
\chi \otimes \chi$. Since
$$(\chi\otimes \chi\otimes \id)(\Phi_0)= (\chi\otimes \id
\otimes \chi)(\Phi_0)= (\id\otimes \chi\otimes \chi)(\Phi_0)=g,$$
application of $\chi\otimes \id\otimes \chi$ gives
\begin{eqnarray*}
\lefteqn{\xi\eta(x)=\xi((\id\ot \chi)(\Delta(x)))}\\
& = & (\chi\ot \id \ot \chi)((\id\ot \Delta)\Delta(x))\\
& = &
(\chi\ot \id \ot \chi)(\Phi_0(\Delta\ot \id)\Delta(x)\Phi_0^{-1})\\
& = & g(\chi\ot \id \ot \chi)(\Delta\ot \id)\Delta(x))
g^{-1}\\
& = & g\eta\xi(x)g^{-1},
\end{eqnarray*}
which is equivalent to $\xi\eta(x)=-\eta\xi(x)$.

Also, application of $\chi\otimes \chi\otimes \id$ yields
\begin{eqnarray*}
\lefteqn{\xi^2(x)=\xi((\chi\ot \id)(\Delta(x)))}\\
& = & (\chi\ot \chi \ot \id)((\id\ot \Delta)\Delta(x))\\
& = & (\chi\ot \chi \ot \id)(\Phi_0(\Delta\ot
\id)\Delta(x)\Phi_0^{-1})=gxg^{-1}=-x,
\end{eqnarray*}
and similarly
application of $\id\otimes \chi\otimes \chi$ yields
$\eta^2(x)=-x$.
\end{proof}

Let $L_g$ be the operator of left multiplication by $g$ in $H_1$.
Then $L_g^2=\id$, and $L_g\xi=-\xi L_g$, $L_g\eta=-\eta L_g$.
Thus, using Lemma \ref{lem1}, we see that the operators
$\xi,\eta,L_g$ define an action on $H_1$ of the Clifford algebra
$Cl_3$ of a $3-$dimensional inner product space.

The following lemma is standard.

\begin{lemma}\label{Cl3}
The algebra $Cl_3$ is semisimple and has two irreducible representations $W_\pm$,
which are both 2-dimensional. They are spanned by
elements $x$ and $gx$ with $\eta(x)=ix$, $\xi(x)=\pm gx$.
\end{lemma}

(Here we abuse the notation by writing $g$ instead of $L_g$.)

Lemma \ref{Cl3} implies that as a $Cl_3$-module, $H_1=nW_+\oplus
mW_-$.

Let $x\in W_\pm\subseteq H_1$ be an eigenvector of $\eta$ with
eigenvalue $i$. Then we have
\begin{eqnarray*}
\lefteqn{\Delta(x)=x\otimes (p_+ \pm ip_-) +p_+\otimes x + p_-\otimes
gx}\\
& = & x\otimes (p_+\pm ip_-) +1\otimes p_+x + g\otimes p_-x.
\end{eqnarray*}

Together with Lemma \ref{sta} (ii), this implies the following
proposition:

\begin{proposition}\label{prop1}
The quasi-Hopf algebra $\gr(H)$ is generated by elements $g$ of
degree $0$, and $x_1,...,x_n,y_1,...,y_m$ of degree $1$, which
satisfy the relations $gx_j=-x_jg, gy_j=-y_jg$, $g^2=1$ (and
possibly other relations). The coproduct of $\gr(H)$ is defined by
the formulas $$ \Delta(g)=g\otimes g, $$ $$ \Delta(x_j)=x_j\otimes
(p_++ip_-) +1\otimes p_+x_j + g\otimes p_-x_j, $$ $$
\Delta(y_j)=y_j\otimes (p_+-ip_-) +1\otimes p_+y_j + g\otimes
p_-y_j. $$
\end{proposition}

\begin{lemma}
In $\gr(H)$, we have $x_j^4=y_l^4=0$ for all $j,l$.
\end{lemma}

\begin{proof}
It suffices to show that $x_j^4=0$, the case of $y_l$ is obtained
by changing $i$ to $-i$. Let $x=x_j$. Using Proposition
\ref{prop1} we find that
$$\Delta(x^2)= x^2\otimes g+ [(1+i)(p_+\otimes p_++p_-\otimes
p_+)+(1-i)(p_+\otimes p_--p_-\otimes p_-)] (x\otimes x)+g\otimes
x^2,$$ hence
$$\Delta(x^4)=x^4\otimes 1+1\otimes x^4.$$
But a finite dimensional quasi-Hopf algebra
cannot contain nonzero primitive elements
(this is proved as in the Hopf case).
Thus, $x^4=0$.
\end{proof}

\begin{lemma}
The numbers $m$ and $n$ are either $0$ or $1$.
\end{lemma}

\begin{proof} Clearly, it suffices to show that $n=0$ or $1$.
Assume the contrary, i.e. that $n\ge 2$.

Introduce the element $z:=g(x_2x_1-ix_1x_2)$. Using Proposition
\ref{prop1} it is checked directly that $$\Delta(z)=z\otimes
1+1\otimes z+ 2(gx_2\otimes p_+x_1+ix_2\otimes p_-x_1). $$ Let $N$
be the smallest integer such that $z^N=0$. It exists since
$\gr(H)$ is finite dimensional. Taking the coproduct of this
equation and looking at the terms of bidegree $(2N-1,1)$, we find
$$ \sum_{k=1}^Nz^{k-1}x_2z^{N-k}=0.
$$ Let us now apply the coproduct to this equality, and look at
the terms of bidegree $(2N-2,1)$ which have a factor $x_2$ in the
second component. This yields $Nz^{N-1}=0$, which is a
contradiction with the minimality of $N$.
\end{proof}

\begin{lemma}
(i) If $m=n=0$, $\gr(H)=H(2)$.

(ii) If $n=1,m=0$ then $\gr(H)=H_+(8)$.

(iii) If $n=0,m=1$ then $\gr(H)=H_-(8)$.

(iv) If $n=1,m=0$ then $\gr(H)=H(32)$.
\end{lemma}

\begin{proof} Statements (i)-(iii) of the lemma
follow from the arguments above, since it is easy to check that
the algebras $H_\pm(8)$ do not have nontrivial graded quasi-Hopf
ideals which do not intersect with $H_1$. It remains to prove
statement (iv).

Assume $m=n=1$, and set $x:=x_1,y:=y_1$. Consider the element
$z:=xy+iyx$. A direct computation shows that $$ \Delta(z)=z\otimes
1+g\otimes z.$$ Since $gz=zg$ and $\gr(H)$ is finite dimensional,
we must have $z=0$. Thus, $\gr(H)$ is a quotient of $H(32)$. But
it is easy to check directly that $H(32)$ does not have nontrivial
graded quasi-Hopf ideals which do not intersect with $H_1$. We are
done.
\end{proof}

Finally, let us show that $H_+(8)$ is not twist equivalent to
$H_-(8)$. Suppose they are. Then there exists a twist $J$ for
$H_-(8)$, and an isomorphism of algebras \linebreak $f: H_+(8)\to
H_-(8)^J$. Such an isomorphism obviously preserves filtration
by powers of the radical, so
we can take its degree zero part. Thus we can assume, without loss
of generality, that $f$ preserves the grading and $J$ is of degree
zero. Then $f$ is the identity on the degree zero part, and
$\Phi^J=\Phi$. So $J$ is a twist of $\Bbb C[\Bbb Z_2]$, hence a
coboundary, and thus we may assume that $J=1\ot 1$ (by changing
$f$). But then $f$ cannot exist, since $S^2(a)=ia$ on $H_+(8)_1$
and $S^2(a)=-ia$ on $H_-(8)_1$.

This completes the proof of Theorem \ref{main3}.

\section{Relation to pointed Hopf algebras}

Let $H=\oplus_{k\ge 0}H_k$ be a graded quasi-Hopf algebra
with radical of codimension 2
and nontrivial associator
on 1-dimensional representations.

The main result of this section is that it is actually possible to
embed $H$ into a twice bigger quasi-Hopf algebra $\tilde H$, which
is twist equivalent to a basic Hopf algebra $H'$ with
$H'/\Rad(H')=\Bbb C[\Bbb Z_4]$. This fact should have
generalizations to the case of general graded basic quasi-Hopf
algebras, which may facilitate applications to the quasi-Hopf case
of known deep results about pointed Hopf algebras.

To construct $\tilde H$, let us adjoin a new element $a$ to $H$
which is grouplike, $a^2=g$, and $aza^{-1}=i^kz$ for $z\in H_k$.
Then the algebra $\tilde H$ generated by $a$ and $H$ (of twice
bigger dimension than that of $H$) is a quasi-Hopf algebra graded
by nonnegative integers, with $\tilde
H_0=\mathbb{C}[\mathbb{Z}_4]$.

\begin{proposition} $\tilde H$ is twist equivalent to a Hopf algebra.
\end{proposition}

\begin{proof}
We claim that the image of the associator $\Phi_0$ in $\tilde
H_0^3$ is homologically trivial. This is because the natural map
$f': H^3(\mathbb{Z}_2,\mathbb{C}^*)\to
H^3(\mathbb{Z}_4,\mathbb{C}^*)$ induced by the projection $f:
\mathbb{Z}_4\to \mathbb{Z}_2$ is zero. Indeed, for a cyclic group
$G$ one has $H^3(G,\mathbb{C}^*)=G^*\otimes G^*$ (functorially in
$G$). Thus for any morphism of cyclic groups $f: G_1\to G_2$, the
induced map of third cohomology groups is $f'=f^*\otimes f^*$. In
our case $f: \mathbb{Z}_4\to \mathbb{Z}_2$, so $f^*:
\mathbb{Z}_2\to \mathbb{Z}_4$ is given by $f^*(1)=2$. Hence,
$(f^*\otimes f^*)(1\otimes 1)=2\otimes 2=1\otimes 4=0$, as
desired.

Thus, by a twist in $\tilde H_0^2$, we can kill $\Phi_0$. So
$\tilde H$ is twist equivalent to a graded pointed
Hopf algebra $A$ generated in degree 1 with $A/Rad(A)=\mathbb{Z}_4$.
\end{proof}

\begin{example} Let $A_\pm(16)$ (for each choice of sign) be
the Hopf algebra of dimension 16 generated by $a,x$ with relations
$a^4=1, axa^{-1}=ix, x^4=0$, such that $a$ is grouplike and
$\Delta(x)=1\otimes x+x\otimes a^{\mp 1}$.

Let $A(64)$ be the Hopf algebra of dimension 64 generated by
$a,x,y$ with relations $a^4=1, axa^{-1}=ix, aya^{-1}=iy, x^4=0,
y^4=0, xy+iyx=0$, such that $a$ is grouplike and
$\Delta(x)=1\otimes x+x\otimes a^{-1}$, $\Delta(y)=1\otimes
y+y\otimes a$.

In other words, $A_\pm(16)=U_q({\frak b}_+)$, where ${\frak b}_+$
is the Borel subalgebra in ${\frak sl}_2$, and $q=\exp(\pm\pi
i/4)$, while $A(64)=\gr(U_q({\frak sl_2})^*)$, for the same $q$.
\end{example}

\begin{proposition}
(i) $\widetilde{H_\pm(8)}$ is twist equivalent to $A_\pm(16)$.

(ii) $\widetilde{H(32)}$ is twist equivalent to $A(64)$.
\end{proposition}

\begin{proof} (i) Consider the Hopf algebra
$A$ obtained by twisting away the associator in
$\widetilde{H_+(8)}$. Then $A=\bigoplus_{k\ge 0} A_k$ is a
$16-$dimensional graded basic (and pointed) Hopf algebra generated
in degree 1. Grouplike elements of this Hopf algebra form a group
$\Bbb Z_4$, and $A_1$ is a free rank 1 module over $A_0$ under
left multiplication. So $A=A_+(16)$ or $A=A_-(16)$. To decide the
sign, assume that $J$ is a pseudotwist for $A_s(16)$, and $f:
\widetilde{H_+(8)}\to A_s(16)^J$ is an isomorphism (where $s$ is a
choice of sign). Since $f$ preserves the filtration by powers of
the radical, we may assume (by taking the degree $0$ part) that
$f$ preserves the grading and $J$ has degree $0$. Then
$\Phi=(1\otimes 1\otimes 1)^J$. This implies that if $Q:=\sum
d_jS(e_j)$, where $J^{-1}=\sum d_j\otimes e_j$, then
$QS(Q)^{-1}=a^2$. Hence the squares of the antipodes in
$\widetilde{H_+(8)}$ and $A_s(16)$ differ by conjugation by $a^2$,
i.e. by a sign on the degree 1 component. For
$\widetilde{H_+(8)}$, the eigenvalue of $S^2$ in degree 1 is $i$.
Thus $s=+$.

(ii) Consider the Hopf algebra $A$ obtained by twisting away the
associator in $\widetilde{H(32)}$. This is a 64-dimensional
graded basic (and pointed) Hopf algebra, which is a quotient of the
amalgamated product of $A_+(16)$ (generated by $a,x$) and
$A_-(16)$ (generated by $a,y$). It is easy to see that $xy+iyx$ is
a primitive element, so it must be zero. Thus, $A=A(64)$, which is
the quotient of this amalgamated product by the relation
$xy+iyx=0$.
\end{proof}


\begin{thebibliography}{AEG}

\bibitem[AS]{as} N. Andruskiewitsch and H-J. Schneider,
Pointed Hopf algebras. New directions in Hopf algebras, 1--68,
{\em Math. Sci. Res. Inst. Publ.}, {\bf 43}, Cambridge Univ.
Press, Cambridge, 2002.

\bibitem[D]{d} V. Drinfeld, Quasi-Hopf algebras. (Russian) {\em Algebra i
Analiz} {\bf 1} (1989), no. 6, 114--148; translation in {\em
Leningrad Math. J.} {\bf 1} (1990), no. 6, 1419--1457.

\bibitem[E]{e} P. Etingof, On Vafa's theorem for tensor categories,
{\em Math. Res. Lett.} {\bf 9} (2002), no. 5-6, 651--657.

\bibitem[ENO]{eno} P. Etingof, D. Nikshych and V. Ostrik, On
fusion categories, {\em submitted}, math.QA/0203060.

\bibitem[EO]{eo} P. Etingof and V. Ostrik, Finite tensor
categories, {\em submitted}, math.QA/0301027.

\bibitem[N]{n} W. Nichols, Bialgebras of type one, {\em Comm. Algebra} {\bf 6}
(1978), no. 15, 1521--1552.

\bibitem[S]{s} M. Sweedler, Hopf Algebras, Benjamin Press, 1968.

\bibitem[TY]{ty} D. Tambara and S. Yamagami,
Tensor categories with fusion rules of self-duality for finite
abelian groups, {\em J. Algebra} {\bf 209} (1998), no. 2,
692--707.

\end{thebibliography}
\end{document}